%

\nopagenumbers
\input amstex
\documentstyle {amsppt}
\magnification=\magstep1
\parindent=20pt
\NoBlackBoxes

\centerline {\bf On the algebra of elementary embeddings of a rank into
itself}
\vskip 1em
\centerline {Richard Laver\footnote{Supported by NSF Grant DMS
9102703.}}
\vskip 2em
\baselineskip=20pt

For $\lambda$ a limit ordinal let ${\Cal E}_\lambda$ Be the set of
elementary embeddings $j : (V_\lambda ,\epsilon ) \rightarrow (V_\lambda
,\epsilon )$, $j$ not the identity.  Then the existence of a $\lambda$
such that ${\Cal E}_\lambda$ is nonempty is a strong large cardinal
axiom (Kanamori--Reinhardt--Solovay [5]).  For $j \in {\Cal E}_\lambda$
let $\text{cr} \,j$ be the critical point of $j$.  Then $\lambda = \sup
\{ j^n (\text{cr} \,j) : n < \omega \}$ (Kunen [7]).

For $j, k \in {\Cal E}_\lambda$ let $j\cdot k = \underset \alpha <
\lambda \to\bigcup j (k \cap V_\alpha )$; then $j\cdot k \in {\Cal
E}_\lambda$ (details are reviewed below).  Writing $jk$ for $j\cdot k$,
we have that $\cdot$ is nonassociative, noncommutative, and satisfies
the left distributive law $a(bc) = (ab)(ac)$.  Letting $j\circ k$ be the
composition of $j$ and $k$, it is seen that $({\Cal E}_\lambda , \cdot ,
\circ )$ satisfies $\sum = \{ a\circ (b\circ c) = (a\circ b)\circ c$,
$(a\circ b)c = a(bc)$, $a(b\circ c) = ab\circ ac$, $a\circ b = ab\circ
a\}$.  $\sum$ implies the left distributive law $(a(bc) = (a\circ b)c =
(ab\circ a)c = (ab)(ac))$.  For $j \in {\Cal E}_\lambda$, let ${\Cal
A}_j$ be the closure of $\{ j\}$ under $\cdot$, and let ${\Cal P}_j$ be
the closure of $\{ j\}$ under $\cdot$ and $\circ$.  A theorem of [9] is 
that ${\Cal A}_j \cong {\Cal A}$, the free
left distributive algebra on one generator, and ${\Cal P}_j \cong {\Cal
P}$, the free algebra on one generator satisfying $\sum$.

For $j \in {\Cal E}_\lambda$ let $\text{cr} \,{\Cal A}_j = \{ \text{cr}
\,k : k \in {\Cal A}_j \}$.  Let $\kappa_0 = \text{cr} \,j$,
$\kappa_{n+1} = j(\kappa_n )$.  Then $\kappa_n = \text{cr} \,(j^n j) \in
\text{cr} \,{\Cal A}_j$.  But there are other members of $\text{cr}
\,{\Cal A}_j$; $\text{cr} \,([(jj)j](jj))$ is between $\kappa_2$ and
$\kappa_3$.  Let $f(n)$ be the cardinality of $\{ \gamma \in \text{cr}
\,{\Cal A}_j : \kappa_n < \gamma < \kappa_{n+1} \}$.  Then $f(0) = 0$,
$f(1) = 0$, $f(2) = 1$, but the author, in calculations, noted that
there are many critical points of members of ${\Cal A}_j$ between 
$\kappa_3$ and $\kappa_4$.

At present an upper bound for $f(3)$ hasn't been computed, but Dougherty
has calculated a very large lower bound.  In addition, Dougherty [4]
proved that $f$ eventually dominates the Ackermann function and is thus
not primitive recursive.  Missing still was a proof that the $f(n)$'s,
even for $n=3$, are finite, that is, that the order type of $\text{cr}
\,{\Cal A}_j$ is $\omega$.  We shall prove this here.  It is a joint
result with Steel.  Namely, Steel in [14] proved that if $k_n \in {\Cal
E}_\lambda$ ($n < \omega$) and $\ell_n = (((k_0 k_1 )k_2 ) \cdots
k_{n-1} )k_n$, then $\sup \{ \text{cr} \ell_n : n < \omega \} = \lambda$.
This is gotten as a simpler special case of the nontrivial studies 
undertaken in
[13], [11], [12], [14].  Steel's result plus results in this paper give
that each $f(n)$ is finite.

The above methods yield a sequence of algebras $\langle P_n : n < \omega
\rangle$, $\text{Card} \,P_n = 2^n$ (with ${\Cal A}_j , {\Cal P}_j$
subalgebras of the inverse limit of $\langle P_n : n < \omega
\rangle$), in which computations of the critical points can be carried
out.  It follows that if $w$ and $u$ are terms in one variable $x$ in 
the language
of $\cdot$ and $\circ$, and $w[j], u[j]$ the corresponding members of
${\Cal P}_j$, then the question, whether or not $\text{cr} \,(w[j]) \leq
\text{cr} \,(u[j])$, is recursive and independent of $j$.

Lastly it is shown that if $k, \ell \in {\Cal A}_j$ and $k \restriction
\text{cr} \,{\Cal A}_j = \ell \restriction \text{cr} \,{\Cal A}_j$, then
$k = \ell$.

The paper is organized as follows.  Lemma 1--Theorem 4 are preliminaries
and a recalling of some results from [9].  In showing ${\Cal A}_j \cong
{\Cal A}$, ${\Cal P}_j \cong {\Cal P}$, a number of combinatorial facts
about ${\Cal A}$ and ${\Cal P}$ were proved.  One of them, the
irreflexivity property, was derived from a large cardinal axiom.
Recently Dehornoy [3] has proved irreflexivity in ZFC; thus
the known results on the structure of ${\Cal A}$ and ${\Cal P}$ are true
in ZFC.  Theorems 5 and 6 summarize these properties.  Lemma 7--Theorem
13 give the abovementioned results about embeddings.

We assume familiarity with large cardinals and elementary embeddings
([5], [6]).  $\lambda$ will always be a limit ordinal.  If $j \in {\Cal
E}_\lambda$, extend $j$ to a map $j : V_{\lambda +1} \rightarrow
V_{\lambda +1}$ by defining, for $A \subseteq V_\lambda$, $j(A) = 
\underset \alpha < \lambda \to\bigcup j (A \cap V_\alpha )$.
A special case is where $A$, as a set of ordered pairs, is a $k \in
V_\lambda$.  That $j(k) \in {\Cal E}_\lambda$ may be checked directly;
it also follows from

\proclaim{Lemma 1}  For $A \subseteq V_\lambda$, $j : (V_\lambda
,\epsilon ,A) \rightarrow (V_\lambda ,\epsilon , jA)$ is elementary.
\endproclaim

\demo{Proof}  Suppose $\Psi$ is a first order sentence in the language
of $(V_\lambda ,\epsilon ,C, c_0 ,\ldots ,c_r )$, where $C \subseteq
V_\lambda$ and each $c_i \in V_\lambda$.  Expanding the language, let
$\Psi^\ast = \forall \,x_0 \cdots \forall \,x_m \Psi '$, $\Psi '$
quantifier free, be the Skolemization of $\Psi : (V_\lambda ,\epsilon
,C, c_0 ,\ldots ,c_r ) \models \Psi$ iff there are $f_i : (V_\lambda
)^{m_i} \rightarrow V_\lambda$ such that $(V_\lambda ,\epsilon , C, c_0
,\ldots ,c_r ,f_0 ,\ldots ,f_n ) \models \Psi^\ast$.  We have that
$(V_\lambda ,\epsilon , C, c_0 ,\ldots ,c_r ,f_0 ,\ldots ,f_n ) \models
\Psi^\ast$ iff, letting ${\Cal T}_0 (x_0 ,\ldots ,x_m ),\ldots ,{\Cal
T}_p (x_0 ,\ldots ,x_m )$ be the terms occurring in $\Psi '$, then for
eventually all/unboundedly many $\delta < \lambda$, $(V_\delta ,\epsilon
,C \cap V_\delta , c_0 ,\ldots ,c_r ,f_0 \cap V_\delta ,\ldots ,f_n \cap
V_\delta ) \models \Psi_\delta^\ast$, where $\Psi_\delta^\ast$ is the
statement $\forall \,x_0 \in V_\delta \cdots \forall \,x_m \in V_\delta
\Big[ \Big( \underset i \leq p \to{\bigwedge\!\!\!\!\!\bigwedge} {\Cal T}_i
(x_0 ,\ldots ,x_m ) ~\text{exists} \Big) \rightarrow \Psi ' (x_0 ,\ldots
,x_m )\Big]$.  To prove the lemma, we are given $(V_\lambda ,\epsilon ,
A, a_0 ,\ldots ,a_r ) \models \Psi$.  Pick $f_0 ,\ldots ,f_n$ such that
$(V_\lambda ,\epsilon ,A, a_0 ,\ldots ,a_r ,f_0 ,\ldots ,f_n ) \models
\Psi^\ast$.  Then for eventually all $\delta < \lambda$, $(V_\delta
,\epsilon , A \cap V_\delta , a_0 ,\ldots ,a_r ,f_0 \cap V_\delta
,\ldots ,f_n \cap V_\delta ) \models \Psi_\delta^\ast$, whence by
elementarity $(V_{j\delta} ,\epsilon , jA \cap V_{j\delta} , ja_0
,\ldots ,ja_r , jf_0 \cap V_{j\delta} ,\ldots ,jf_n \cap V_{j\delta} )
\models \Psi_{j\delta}^\ast$.  Thus $(V_\lambda ,\epsilon ,jA, ja_0
,\ldots ,ja_r , jf_0 ,\ldots ,jf_n ) \models \Psi^\ast$, so $(V_\lambda
,\epsilon ,jA, ja_0 ,\ldots ,ja_r ) \models \Psi$.

For $j, k \in {\Cal E}_\lambda$ let $j\cdot k$ (abbreviated $jk$) $=
j(k)$.  Let $j\circ k$ be the composition of $j$ and $k$.  For $k_0 ,k_1
,\ldots ,k_n \in {\Cal E}_\lambda$ let $k_0 k_1 k_2 \cdots k_n =
(((k_0 k_1 )k_2 )\cdots )k_n$ and let $k_0 k_1 k_2 \cdots k_{n-1}
\circ k_n = ((((k_0 k_1 )k_2 )\cdots )k_{n-1} )\circ k_n$.  Write $k =
k_0 k_1 \cdots k_{n-1} \ast k_n$ to mean that $k = k_0 k_1 \cdots
k_{n-1} k_n$ or $k = k_0 k_1 \cdots k_{n-1} \circ k_n$.  Make these
conventions also for other algebras on operations $\cdot$ and $\circ$.
For $\theta < \lambda$ let $\log_k (\theta )$ be the least $\mu$ with
$k(\mu ) \geq \theta$.  Define $k \overset \ast \to\bigcap V_\theta$ to
be $\{ \langle x,y \rangle \in V_\theta \times V_\theta , y \in k(x)\}$.
Let $k \overset \theta \to{=} \ell$ mean that $k \overset \ast
\to\bigcap V_\theta = \ell \overset \ast \to\bigcap V_\theta$, and let
$k \underset \theta \to{=} \ell$ if $k \restriction V_\theta = \ell
\restriction V_\theta$.
\enddemo

The following is a list of properties of ${\Cal E}_\lambda$ which will
be used without comment below.  If $k \in {\Cal E}_\lambda$, 
$\alpha \geq \text{cr} \,k$,
then $k(\alpha ) > \alpha$.  For $k,j \in {\Cal E}_\lambda$, $\text{cr}
\,kj = k(\text{cr} \,j)$.  If $\alpha < \text{cr} \,k$, then $kj(\alpha
) = k(j\alpha )$.  Define $j(< d) = \sup  \{ j\beta : \beta < \alpha
\}$.  We have $j(< \alpha) < j(\alpha )$ iff $\text{cf}
\,\alpha \geq \text{cr} \,j$; $(j\circ k)(< \alpha ) = j(< k(< \alpha
))$.  If $k, \ell_0 ,\ell_1 ,\ldots ,\ell_n \in {\Cal E}_\lambda$ and
$\text{cr} \,k = \gamma$, then $k\ell_0 \ell_1 \cdots \ell_{n-1} \ast
\ell_n \overset \gamma \to{=} \ell_0 \ell_1 \cdots \ell_{n-1} \ast
\ell_n$.

Let $\mu = \log_k \theta$.  Then if $k \overset \theta \to{=} \ell$,
then $\mu = \log_\ell \theta$.  $k \underset \mu +1 \to{=} \ell
\Rightarrow k \overset \theta \to{=} \ell \Rightarrow k \underset \mu
\to{=} \ell$.  $k\ell \overset \ast \to\bigcap V_\theta = k(\ell
\overset \ast \to\bigcap V_\mu ) \overset \ast \to\bigcap V_\theta$.
If $k \overset \theta \to{=} k'$, $\ell \overset \theta \to{=} \ell '$,
then $k\circ \ell \overset \theta \to{=} k' \circ \ell '$.  If $k
\overset \theta \to{=} k'$, $\ell \overset \mu \to{=} \ell '$, then
$k\ell \overset \theta \to{=} k' \ell '$.  

\proclaim{Lemma 2 ([9])}  Suppose $k, \ell_0 ,\ell_1 ,\ldots ,\ell_n \in
{\Cal E}_\lambda$.  Let $\text{cr} \,k = \gamma$ and for $m \leq n$,
$\theta_m = \min \{ k\ell_0 \cdots \ell_i (\gamma ) : i < m\}$ (say
$\theta_m = \lambda$ if $m = 0$).  Then
\item {(i)} $k\ell_0 \cdots \ell_n \overset \theta_n \to{=} k(\ell_0
\cdots \ell_n )$.
\item {(ii)} $\theta_n = \min \{ k(\ell_0 \cdots \ell_i )(\gamma ) : i <
n\}$.
\item {(iii)} if for all $m < n$, $\text{cr} \,(\ell_0 \cdots \ell_m )
< \gamma$, then
\itemitem {(a)} $\theta_n > k(\gamma )$.
\itemitem {(b)} if $\delta < \gamma$ and $(\ell_0 \cdots \ell_n )(\delta
) = \gamma$, then $k\ell_0 \cdots \ell_n (\delta ) = k(\ell_0 \cdots
\ell_n )(\delta ) = k(\gamma )$.
\endproclaim

\demo{Proof}  (i) By induction on $n$   The case $n=0$ is by definition.
Suppose it holds for $n-1$.  Since $\text{cr} \,k = \gamma$, $\ell_n
\overset \gamma \to{=} k\ell_n$.  Applying $k\ell_0 \cdots \ell_{n-1}$
to both sides and using the induction hypothesis $k\ell_0 \cdots
\ell_{n-1} \ell_n \overset k\ell_0 \cdots \ell_{n-1} (\gamma ) \to{=}
k\ell_0 \cdots \ell_{n-1} (k\ell_n ) \overset \theta_{n-1} \to{=}
k(\ell_0 \cdots \ell_{n-1} )(k\ell_n ) = k(\ell_0 \cdots \ell_{n-1}
\ell_n )$.  Since $\theta_n = \min \{ k\ell_0 \cdots \ell_{n-1} (\gamma
), \theta_{n-1} \}$, we are done.
\enddemo

(ii) By induction on $n$, using (i).  If either of $k\ell_0 \cdots
\ell_{n-1} (\gamma )$, $k(\ell_0 \cdots \ell_{n-1} )(\gamma )$ is less
than $\theta_{n-1}$, then $k\ell_0 \cdots \ell_{n-1} (\gamma ) =
k(\ell_0 \cdots \ell_{n-1} )(\gamma ) = \theta_n$.  Otherwise, $\min \{
k(\ell_0 \cdots \ell_i )(\gamma ) : i < n\} = \min \{ k\ell_0 \cdots
\ell_i (\gamma ) : i < n\} = \theta_{n-1} = \theta_n$.

(iii) For (a), $\theta_n = k(\ell_0 \cdots \ell_i )(\gamma )$ for some
$i < n$.  And $k(\ell_0 \cdots \ell_i )(\gamma ) >$\newline
$k(\ell_0 \cdots
\ell_i )[(\log_{\ell_0 \cdots \ell_i} \gamma ) +1] = k(\ell_0 \cdots
\ell_i )(k[(\log_{\ell_0 \cdots \ell_i} \gamma ) +1]) = k[\ell_0 \cdots
\ell_i [(\log_{\ell_0 \cdots \ell_i} \gamma ) +1]] > k(\gamma )$.  Part
(b) follows form (i) and (iii)(a).

\proclaim{Corollary 3 ([9])}  If $n > 0$, $k_0 ,k_1 ,\ldots ,k_n \in
{\Cal E}_\lambda$, then $k_0 \not= k_0 k_1 \cdots k_{n-1} \ast k_n$.
\endproclaim

\demo{Proof}  Assume $k_0 = k_0 k_1 \cdots k_n$.  Let $\gamma = \max \{
\text{cr} \,(k_0 k_1 \cdots k_i ) : 0 \leq i < n\}$, and let $\ell$ be a
$k_0 k_1 \cdots k_i$ with $\text{cr} \,\ell = \gamma$.  Since $\ell
k_{i+1} \cdots k_n k_0 \cdots k_i = \ell$, there would be $u_0 ,u_1
,\ldots ,u_r$ with $\text{cr} \,\ell = \gamma$, $\text{cr} \,(\ell u_0
\cdots u_t ) < \gamma$ ($t < r$), and $\text{cr} \,(\ell u_0 \cdots u_r
) = \gamma$.  But this is impossible by Lemma 1.  I.e., $\gamma =
\text{cr} \,(\ell u_0 \cdots u_r ) = \ell u_0 \cdots u_{r-1} (\text{cr}
\,u_r )$, so $\text{cr} \,u_r < \gamma$ since $\text{cr} \,(\ell u_0
\cdots u_{r-1} ) < \gamma$.  Since $\theta_{r-1} > \gamma$, we have $\ell 
u_0 \cdots u_{r-1} (\text{cr} \,u_r ) = \ell (u_0 \cdots u_{r-1} )(\ell
(\text{cr} \,r)) = \ell (u_0 \cdots u_{r-1} (\text{cr} \,u_r ))$, but
$\gamma \notin \text{range} \,\ell$, a contradiction.  And if $k_0 = k_0
k_1 \cdots k_{n-1} \circ k_n$, then $k_0 k_1 = k_0 k_1 \cdots k_{n-1}
(k_n k_i )$, which reduces to the previous case.

Recall that ${\Cal A}$ is the free algebra on one generator $x$ satisfying 
the left distributive law, and ${\Cal P}$ is the free algebra on one
generator $x$ satisfying $\sum$.  Then ${\Cal A} = A/\equiv$, where $A$
is the set of terms in the language of $\cdot$ in one variable $x$, and
for $w, v \in A$, $w \equiv v$ if $w$ and $v$ are equivalent using the
left distributive law.  Similarly ${\Cal P} = P/\equiv$, where $P$ is
the set of terms in $x$ in the language of $\cdot$ and $\circ$, and
$\equiv$ means equivalent using $\sum$.  The same notation $\equiv$ is
used because if $u, v \in A$ and $u \equiv v$ via $\sum$, then $u \equiv
v$ via the left distributive law ([9], section 1).  So ${\Cal A}$
may be viewed as a subalgebra of ${\Cal P}$ restricted to $\cdot$ .

\proclaim{Theorem 4 ([9])}  For $j \in {\Cal E}_\lambda$, ${\Cal A}_j 
\cong {\Cal A}$, ${\Cal P}_j \cong {\Cal P}$.
\endproclaim

The only application of this below is via the linear ordering of Theorem
5, but we recount some general facts about ${\Cal A}$ and ${\Cal P}$.
Every member $w$ of ${\Cal P}$ is a composition $w_1 \circ w_2 \circ
\cdots \circ w_{c(w)}$, each $w_i \in {\Cal A}$, where $c(w)$ only
depends on $w$.  If $w \in {\Cal P}$, let $w^{(0)} = w$, $w^{(i+1))} =
w^{(i)} w^{(i)}$.  Then for $i \leq n$, $w^{(i)} w^{(n)} = w^{(n+1)}$ by
induction using the left distributive law.

For $u, w \in {\Cal P}$ write $u < w$ if $w$ can be written as a term of
length greater than one in members of ${\Cal P}$, at least one of which
is $u$.  Write $u <_L w$ if $u$ occurs on the left of such a term:  $w =
uu_0 u_1 \cdots u_{n-1} \ast u_n$.  Then $<_L$ is transitive and $x <_L
w$ for all $w \in {\Cal P} - \{ x\}$.  If $w <_L v$, then $aw <_L a
\circ w <_L av$.  If $p \in {\Cal P}$, then it is seen by induction that
for sufficiently large $n$, $p <_L x^{(n)}$ and $px^{(n)} =
x^{(n+c(p))}$.

In [9] it was shown that $<_L$ is a linear ordering on ${\Cal A}$ and
${\Cal P}$.  The statement that $<_L$ is irreflexive was derived as a
consequence of a large cardinal axiom:  ~it follows from Corollary 3.
The remaining part of the linear ordering theorem involved a normal form
for the members of ${\Cal A}$ and ${\Cal P}$.  

Dehornoy ([1], [2]) independently proved the remaining part of the
linear ordering theorem (that for all $a, b \in A$ at least one of $a
<_L b$, $a = b$, $b <_L a$ holds) in ZFC, by a different method.  For
$u, v \in A$, let $u \rightarrow v$ mean that there is a sequence $u =
u_0 ,u_1 ,\ldots ,u_n = v$ where $u_{i+1}$ comes from $u_i$ by replacing
a component of the form $a(bc)$ by $(ab)(ac)$.  Then ([1]) the rewriting
rules for ${\Cal A}$ are confluent:  ~if $u, v \in A$, $u \equiv v$,
then for some $r \in A$, $u \rightarrow r$ and $v \rightarrow r$.  More
recently ([3]) Dehornoy has proved the irreflexivity of $<_L$ in ZFC and
found a representation of ${\Cal A}$ in terms of an operation on certain
braid words. Theorem 5
contains results about ${\Cal A}$ and ${\Cal P}$ from ([9]), which 
used irreflexivity and which are now known from ZFC via 
Theorem 6.

For $u \in {\Cal P}$ a $u$-prenormal sequence is a term $u_0 u_1 \cdots
u_{n-1} \ast u_n$ ($n \geq 0$) of members of ${\Cal P}$, where $u=u_0$
and $u_{i+2}
\leq_L u_0 u_1 \cdots u_i$ for all $i \leq n-2$, and in the case $\ast =
\circ$, $u_n <_L u_0 u_1 \cdots u_{n-2}$.

\proclaim{Theorem 5} On ${\Cal A}$ and ${\Cal P}$,
\item {(i)} $<_L$ is a linear ordering.
\item {(ii)} $uv = uw$ iff $v = w$, $uv <_L uw$ iff $v <_L w$.
\item {(iii)} the word problems are decidable.
\item {(iv)} if $u \leq_L v$, then there is a (unique) $u$-prenormal
sequence equal to $v$.  Moreover, let $(T,\ell )$ be the rooted,
labelled tree such that $\ell (\text{root} \,T) = v$, $\ell (t) \leq_L u
\Rightarrow t$ is maximal in $T$, $u <_L \ell (t) = u_0 u_1 \cdots
u_{n-1} \ast u_n$ ($u_0 u_1 \cdots u_{n-1} \ast u_n$ $u$-prenormal)
$\Rightarrow t$'s immediate successors in $T$ are $t_0 ,\ldots ,t_n$
with $\ell (t_i ) = u_i$.  Then $T$ is finite.
\item {(v)} $<_L = <$.
\endproclaim

\noindent {\bf Remarks.}  The linear orderings of ${\Cal A}$ and ${\Cal
P}$ both have order type $\omega \cdot (1+\eta )$.  (ii) and (iii) are
immediate consequences of (i).  (iv) is the division form version [10] of
the normal forms of [9].  (v) was derived by the author from (iv);
McKenzie derived (v) from (i).  To $<_L$ compare $v$ and $w$ in ${\Cal P}$
using (iv), put $v, w$ in $x$-division form (the decomposition given by
(iv) where $u=x$) and compare the two forms lexicographically (see ([9],
[10]) for details).

Let $B_\infty$ be the infinite dimensional braid group (given by
generators $\sigma_1 ,\sigma_2 ,\ldots$ and the relations
$\sigma_n \sigma_i = \sigma_i \sigma_n$ ($|i-n| > 1$), $\sigma_n
\sigma_{n+1} \sigma_n = \sigma_{n+1} \sigma_n \sigma_{n+1}$).  For
$\alpha \in B_\infty$ let $s (\alpha )$ be the image of
$\alpha$ under the shift map $\sigma_i \rightarrow \sigma_{i+1}$.  
Dehornoy's bracket operation [\quad ] on $B_\infty$ is defined by 
$\alpha [\beta ] = \alpha
\cdot s(\beta )\cdot \sigma_1 \cdot (s \alpha )^{-1}$.  This operation
is left distributive.  If ${\Cal B}$ is any left distributive algebra on
one generator $x$ which satisfies left cancellation (as ${\Cal A}$ turns
out to be, because of the linearity of $<_L$),  then the action of a braid
generator on a member of $({\Cal B} )^\omega$, given by $\langle b_1
,b_2 ,\ldots ,b_i ,b_{i+1} ,\ldots \rangle^{\sigma_i} = \langle b_1 ,b_2
,\ldots ,b_i b_{i+1} ,b_i ,\ldots \rangle$, extends to a partial action
of $B_\infty$ on $({\Cal B} )^\omega$.  And for each $b \in {\Cal B}$
there is an $\alpha_b \in B_\infty$ such that $\langle x, x, x,\ldots
\rangle^{\alpha_b} = \langle b, x, x,\ldots \rangle$, given by $\alpha_x
= \varepsilon$, $\alpha_{bc} = \alpha_b [\alpha_c ]$.  Letting 
$\text{cl} \,(\alpha )$ be the closure of $\{ \alpha \}$ under
[\quad ], we thus have a homomorphism from ${\Cal B}$ to $\text{cl}
\,(\varepsilon )$.

\proclaim{Theorem 6 (Dehornoy, [3])}  
\item {(i)} $<_L$ is irreflexive in ZFC.
\item {(ii)} if $\alpha \in B_\infty$, then $(\text{cl} \,(\alpha ),
[\quad ]) \cong {\Cal A}$.
\endproclaim

Also included in [3] are applications to $B_\infty$ itself.  In Larue
[8] a direct proof of one of them is given, which yields a shorter
proof of Theorem 6.

Suppose $j, k \in {\Cal E}_\lambda$.  Define the $n$th iterate $I_n
(j,k)$ of $\langle j,k\rangle$ by (writing $I_n = I_n
(j,k)$) $I_0 = k$, $I_1 = j$, $I_{n+2} = I_{n+1} I_n$.  Then $I_{n+1}
\circ I_n = j\circ k$ by applications of the $a\circ b = ab\circ a$ law.

Let $\text{id}$ be the identity function on $V_\lambda$.  Topologize
$(V_\lambda )^{V_\lambda}$ by taking basic open sets of the form $\{ f
\in (V_\lambda )^{V_\lambda} : f \restriction V_\theta = g\}$ ($\theta <
\lambda$, $g : V_\theta \rightarrow V_\lambda$).  Then ${\Cal E}_\lambda
\cup \{ \text{id} \}$ is a closed subspace of $(V_\lambda
)^{V_\lambda}$.  Namely suppose $k \in (V_\lambda )^{V_\lambda}$,
$\theta_0 < \cdots < \theta_n \cdots ,\sup \{ \theta_n : n < \omega \} =
\lambda$, $k_n \in {\Cal E}_\lambda$ and $k \underset \theta_n \to{=}
k_n$.  Assume without loss of generality $k_0 \not= \text{id}$,
$\text{cr} \,k_0 = \kappa_0 < \theta_0$.  Then $\kappa_{n+1} =
k(\kappa_n )$ is well defined for all $n$, and by Kunen's theorem for
the $k_n$'s, $\sup \{ \kappa_n : n < \omega \} = \lambda$.  Since $k
\restriction V_{\kappa_n} : (V_{\kappa_n} ,\epsilon ) \rightarrow
(V_{\kappa_{n+1}} ,\epsilon )$ is elementary for all $n$, $k \in {\Cal
E}_\lambda$.

A particular case is when $k$ is a direct limit of members of ${\Cal
E}_\lambda$.  That is, if $\ell_n \in {\Cal E}_\lambda$ ($n < \omega$),
let $\underset \longrightarrow \to\lim ~\langle \ell_n : n < \omega
\rangle$ be the direct limit map.  Say that $\underset \longrightarrow
\to\lim ~\langle \ell_n : n < \omega \rangle$ stabilizes if for each $x
\in V_\lambda$ the sequence $\langle (\ell_n \circ \cdots \circ \ell_1
\circ \ell_0 )(x) : n < \omega \rangle$ is eventually constant.  So if
$\underset \longrightarrow \to\lim ~\langle \ell_n : n < \omega \rangle$
stabilizes, then $\underset \longrightarrow \to\lim ~\langle \ell_n : n
< \omega \rangle \in {\Cal E}_\lambda$.

For ${\Cal H} \subseteq {\Cal E}_\lambda$, let ${\Cal A}_{\Cal H}$ be
the closure of ${\Cal H}$ under $\cdot$ , and let ${\Cal P}_{\Cal H}$ be
the closure of ${\Cal H}$ under $\cdot$ and $\circ$.  Let $\overline
{\Cal H}$ be the topological closure of ${\Cal H}$; $\overline {\Cal H}
= \{ k : \forall \,\theta < \lambda ~\exists \,h \in {\Cal H} ~k \underset
\theta \to{=} h \}$.

\proclaim{Lemma 7}  If $j, k \in {\Cal E}_\lambda$, $\text{cr} \,j >
\text{cr} \,k$, $I_n = I_n (j,k)$, then
\item {(i)} $\text{cr} \,I_0 = \text{cr} \,I_2 = \cdots , ~\text{cr}
\,I_1 < \text{cr} \,I_3 < \cdots$
\item {(ii)} $I_{2n+2} \underset \text{cr} \,I_{2n+1} \to{=} j \circ k$
\item {(iii)} $n,m$ odd $\Rightarrow I_n (\text{cr} \,I_m ) < \text{cr}
I_{m+2}$ (whence $\underset i \to\sup (\text{cr} \,I_{2i+1} ) =
\lambda$).
\item {(iv)} if ${\Cal H} \subseteq {\Cal E}_\lambda$, then $\{ p
\restriction V_\theta : \theta < \lambda , ~p \in {\Cal P}_{\Cal H} \} =
\{ p \restriction V_\theta : \theta < \lambda , ~p \in {\Cal A}_{\Cal H}
\}$, whence $\overline {\Cal A}_{\Cal H} = \overline {\Cal P}_{\Cal H}$.
\endproclaim

\demo{Proof}  (i) is clear.  For (ii), if $\Vert x\Vert < \text{cr}
\,I_{2n+1}$, then $I_{2n+2} (x) = I_{2n+1} I_{2n} (x) =$\hfill\break
$I_{2n+1} I_{2n} (I_{2n+1} x)
= I_{2n+1} (I_{2n} x) = I_1 (I_0 x)$ (as $I_{m+1} \circ I_m
= I_1 \circ I_0$).

For (iii), pick $\alpha < \text{cr} \,I_m$ with $I_{n-1} (\alpha ) >
\text{cr} \,I_m$ (possible since $\text{cr} \,I_{n-1} = \text{cr} \,I_0
< \text{cr} \,I_m$).  Then $\text{cr} \,I_{m+2} = I_{m+1} (\text{cr}
\,I_m ) > I_{m+1} (\alpha ) = j\circ k (\alpha )$ by (ii), and $j\circ
k(\alpha ) = I_n \circ I_{n-1} (\alpha ) \geq I_n (\text{cr} \,I_m )$.
In particular, then, $I_1 (\text{cr} \,I_{2n+1} ) < \text{cr}
\,I_{2n+3}$ for all $n$, so $\sup \langle I_{2n+1} : n < \omega \rangle
= \lambda$.

For (iv), let $j, k \in {\Cal E}_\lambda$.  Then there is a $p \in {\Cal
A}_{\{ j,k\}}$ with $p \overset \theta \to{=} (j\circ k)$.  Namely, by
rewriting $j\circ k$ as $jk\circ j$ if necessary,  we may assume
$\text{cr} \,j > \text{cr} \,k$.  Let $I_m = I_m (j,k)$, and pick $m$
with $\text{cr} \,I_{2m+1} \geq \theta$.  Then $I_{2m+2} \underset
\theta \to{=} (j\circ k)$ (so $\overset \theta \to{=} (j\circ k)$) and
$I_{2m+2} = I_{2m+1} (I_{2m-1} (\cdots (I_3 (I_1 I_0 )))) \in {\Cal
A}_{\{ j,k\}}$.  Then for $p \in {\Cal P}_{\Cal H}$ and $\gamma <
\lambda$ there is $q \in {\Cal A}_{\Cal H}$ with $q \underset \gamma
\to{=} p$, by induction on ${\Cal P}_{\Cal H}$ taking $\theta = p(\gamma
)$.  This gives (iv).
\enddemo

We show, for a subsequent paper concerning limits, that if
$\text{cr} \,j > \text{cr} \,k$, $I_n = I_n (j,k)$, then $\underset
\longrightarrow \to\lim ~\langle I_{2n+1} : n < \omega \rangle$
stabilizes.

\proclaim{Lemma 8}  If $j, k \in {\Cal E}_\lambda$ and $\text{cr} \,j >
\text{cr} \,k$ (respectively, $\text{cr} \,j \leq \text{cr} \,k$), then
there is an $s \in {\Cal E}_\lambda$, $s$ a direct limit of members of
${\Cal A}_{\{ j,k\}}$, with $sk = j\circ k$ (respectively, $sj = j\circ
k$).
\endproclaim

\demo{Proof}  By rewriting $j\circ k$ as $jk \circ j$ if $\text{cr} \,j
\leq \text{cr} \,k$, the lemma reduces to the $\text{cr} \,j > \text{cr}
\,k$ case.  Let $I_n = I_n (j,k)$.  By Lemma 7 (ii) and (iii),
$(I_{2n+1} \circ \cdots \circ I_3 \circ I_1 )(k) = I_{2n+2} \underset
\text{cr} \,I_{2n+1} \to{=} j\circ k$, and $\sup \{ \text{cr} \,I_{2n+1}
: n < \omega \} = \lambda$, so it remains to show that $\underset
\longrightarrow \to\lim ~\langle I_{2n+1} : n < \omega \rangle$
stabilizes.  If not, then by Lemma 7 (i) there would be an 
$\alpha < \lambda$ with
$$(I_{2n-1} \circ \cdots \circ I_3 \circ I_1 )(\alpha ) < (I_{2n+1}
\circ I_{2n-1} \circ \cdots \circ I_3 \circ I_1 )(\alpha )$$
for all $n$.  Then by Lemma 7 (iii),
$\sup \{ (I_{2n+1} \circ \cdots \circ I_3 \circ I_1 )(\alpha )
: n < \omega \} = \lambda$.  Pick $N$ such that $k^N (\text{cr} \,k) >
\alpha$.  By Lemma 7 (iii), pick $m$ such that $(j\circ k)^N (\text{cr} 
\,(j\circ k)) <
\text{cr} \,I_{2m+1}$.  Let $\alpha_m = (I_{2m+1} \circ \cdots \circ I_3
\circ I_1 )(\alpha ) \geq \text{cr} \,I_{2m+1}$.  Then $(I_{2m+1} \circ
\cdots \circ I_3 \circ I_1 )(k \cap V_\alpha ) = (I_{2m+1} \circ \cdots
\circ I_3 \circ I_1 )(k) \cap V_{\alpha_m} \supseteq (j\circ k) \cap
V_{\text{cr} \,I_{2m+1}}$.  But the right hand side (unlike the left
hand side) is a map $h$ such that $h^N (\text{cr} \,h)$ exists, a
contradiction.

The inequalities showing that $\underset \longrightarrow \to\lim
\langle I_{2n+1} : n < \omega \rangle$ stabilizes, where $I_m = I_m 
(j,k)$, $\text{cr} \,j > \text{cr} \,k$, can also be directly computed.
For example, let $\theta_n = (I_{2n+1} )^{n+1} (\text{cr} \,I_{2n+1} )$.
Then

\noindent (1) \quad $\theta_n = (j\circ k)(\theta_{n-1} ) < \text{cr} 
\,I_{2n+3}$
\enddemo

\demo{Proof}  For $n = 0$, $\theta_0 = I_1 (\text{cr} \,I_1 ) <
\text{cr} \,I_3$ by Lemma 7 (iii).  Assume (1) is true for $n-1$.  Then
$\text{cr} \,I_{2n+3} = I_{2n+2} (\text{cr} \,I_{2n+1} ) > I_{2n+2}
(\theta_{n-1} ) =(j\circ k)(\theta_{n-1} )$ (by the induction hypothesis
and Lemma 5 (ii))
$$\aligned &= I_{2n}  (I_{2n-1} (\theta_{n-1} )) = I_{2n} (I_{2n-1}
((I_{2n-1} )^n (\text{cr} \,I_{2n-1} ))) \\
&= I_{2n} ((I_{2n-1} )^{n+1} (\text{cr} \,I_{2n-1} )) = (I_{2n+1} )^{n+1}
(\text{cr} \,I_{2n+1} ) = \theta_n \endaligned$$

\noindent (2) \quad $\theta_n = I_{2n+1} (I_{2n-1} (\cdots (I_3 ((I_2 )^n
(\theta_0 )))))$
\enddemo

\demo{Proof}  For $n=0$ the expressions are the same.  Assuming it is
true for $n-1$ we have $\theta_n = (j\circ k)(\theta_{n-1} ) = I_{2n+1}
(I_{2n} (\theta_{n-1}  )) = I_{2n+1} (I_{2n} (I_{2n-1} (I_{2n-3} (\cdots
(I_3 (I_2^{n-1} (\theta_0 )))))))$, which, since $I_{2k} \circ I_{2k-1}
= I_{2k-1} \circ I_{2k-2}$, equals the right hand side of (2).
\enddemo

To see that $\underset \longrightarrow \to\lim ~\langle I_{2n+1} : n <
\omega \rangle$ stabilizes on a given $\alpha < \lambda$, pick $N$ such
that $(I_2 )^N (\theta_0 ) > I_1 (\alpha )$ (possible since $\text{cr}
\,I_2 < \text{cr} \,I_1 < I_1 (\text{cr} \,I_1 ) = \theta_0$).  Then
by (1) and (2),
$$(I_{2N+1} \circ \cdots \circ I_3 \circ I_1 )(\alpha ) < (I_{2N+1}
\circ \cdots \circ I_3 )(I_2^N (\theta_0 )) = \theta_N < \text{cr}
\,I_{2N+3} .$$

\proclaim{Theorem 9 (Steel [14])}  If $k_n \in {\Cal E}_\lambda$ ($n <
\omega$), then $\sup \{ \text{cr} \,(k_0 k_1 \cdots k_n ) : n < \omega
\} = \lambda$.
\endproclaim

\demo{Proof}  For proofs of the basic facts about extenders sketched
here, see Martin and Steel ([11], [12], [14]).  Steel argues as in
Mitchell's proof [13] that $\triangleleft$ is well founded on normal
ultrafilters, except using extenders.  Suppose $\theta < \lambda$.  For
$k \in {\Cal E}_\lambda$ let $E_k^\theta$ be the extender of length
$\theta$ induced by $k$.  That is, for $a \in [\theta ]^{<\aleph_0}$ let
$\mu_a^k = \log_k (\sup a+1)$.  Let $E_a = \{ X \subseteq [\mu_a^k 
]^{\overset = \to{a}} : a \in k(X)\}$, and define $E_k^\theta = \langle
E_a : a \in [\theta ]^{< \aleph_0} \rangle$.  Let $\text{Ult} \,(V, 
E_k^\theta )$
be the ultraproduct of $V$ by $E_k^\theta$.  The elements of $\text{Ult}
\,(V, E_k^\theta )$ are equivalence classes $[f, a]$, where $a \in
[\theta ]^{< \aleph_0}$ and $f : [\mu_a^k ]^{\overset = \to{a}}
\rightarrow V$.  To define when $[f, a] \sim [g,b]$, let $a\cup b$ in
increasing order be $\{ \alpha_0 ,\alpha_1 ,\ldots ,\alpha_r \}$, $a =
\{ \alpha_{i_0} ,\ldots ,\alpha_{i_n} \}$, $b = \{ \alpha_{k_0} ,\ldots
,\alpha_{k_\ell} \}$; then $[f,a] \sim [g,b]$ iff $\{ u : \bar f (u) =
\bar g (u)\} \in E_{a\cup b}$, where $\bar f (\{ \gamma_0 ,\ldots
,\gamma_r \} ) = \{ \gamma_{i_0} ,\ldots ,\gamma_{i_n} \}$, $\bar g (\{
\gamma_0 ,\ldots ,\gamma_r \} ) = \{ \gamma_{k_0} ,\ldots
,\gamma_{k_\ell} \}$.  Then $\text{Ult} \,(V, E_k^\theta )$ is well
founded.  Let $i_k^\theta : V \rightarrow \text{Ult} \,(V, E_k^\theta )$
be the canonical elementary embedding.

If $\theta < \lambda$ is inaccessible, then claim $V_\theta \subseteq
\text{Ult} \,(V, E_k^\theta )$, in fact for each $[f,a]$ as above with
range $f \subseteq V_\theta$, $[f,a] = kf(a)$.  Namely, if $\sigma
,\delta < \theta$, $f : [\sigma ]^n \rightarrow V_\lambda$, $g : [\delta
]^m \rightarrow V_\lambda$, then $[f,a] \in [g,b]$ iff $k\bar f (a\cup
b) \in k \bar g (a\cup b)$.  Then by induction on $\rho < \theta$ one
establishes $V_\rho \subseteq \text{Ult} \,(V, E_k^\theta )$ 
by showing each $kh (\langle \rho ,\alpha \rangle )$ ($\alpha <
\overset = \to{V_\rho} )$ is in $\text{Ult} \,(V, E_k^\theta )$, where
$\langle h\langle\delta ,\alpha\rangle : \alpha < \overset = \to{V_\delta}
\rangle$ is a well ordering of $V_\delta$ ($\delta \leq \log_k \rho$).

It follows that $i_\ell^\theta \overset \theta \to{=} k$ and that if 
$\ell \in {\Cal E}_\lambda$, then $E_{k\ell}^\theta = i_\ell^\theta 
(E_\ell^{\log_k
\theta} ) \cap V_\theta$.  So $E_{k\ell}^\theta \in \text{Ult} \,(V,
E_k^\theta )$ and is, there, the extender of length $\theta$ induced by
$i_\ell^\theta (\ell )$.
Suppose $\theta < \lambda$ is inaccessible, $k, k\ell \in {\Cal
E}_\lambda$ and $\text{cr} \,k < \theta$.  Then we claim that $i_\ell^\theta
(\theta) 
> i_{k\ell}^\theta (\theta )$.  In $\text{Ult} (V, E_k^\theta )$ let 
$i_{k\ell} ' :
\text{Ult} \,(V, E_k^\theta ) \rightarrow \text{Ult} \,(\text{Ult} \,(V,
E_k^\theta ), E_{k\ell}^\theta )$ be the canonical elementary embedding.
Since $\theta$ is inaccessible, for each $a \in [\theta ]^{< \aleph_0}$
every $f : [\mu_a^{k\ell} ] \rightarrow \theta$ lies in $\text{Ult}
\,(V, E_k^\theta )$; thus $i_{k\ell} (\theta ) = i_{k\ell} ' (\theta )$.
Since $\text{cr} \,k < \theta$, $\text{Ult}
\,(V, E_k^\theta ) \models i_\ell^\theta (\theta ) > \theta$ is inaccessible.
Thus $i_\ell^\theta (\theta ) > i_{k\ell} ' (\theta )$.  This proves the claim.

If the theorem failed, pick an inaccessible $\theta < \lambda$ with each
$\text{cr} \,(k_0 \cdots k_n ) < \theta$.  Then by the claim, for each $n$,
$i_{k_0 \cdots k_n}^\theta (\theta) > i_{k_0 \cdots k_{n+1}}^\theta
(\theta )$, a contradiction.

If ${\Cal H} \subseteq {\Cal E}_\lambda$, let $\text{cr} \,{\Cal H} = \{
\text{cr} \,k : k \in {\Cal H} \}$.  Note that since every member of 
${\Cal P}_{\Cal H}$ is a composition of members of ${\Cal A}_{\Cal H}$,
$\text{cr} \,{\Cal P}_{\Cal H} = \text{cr} {\Cal A}_{\Cal H}$ (that fact
is also a consequence of Lemma 5 (iv)).  For $X \subseteq \text{Ord}$,
let $\text{tp} \,X$ be the order type of $X$.

\proclaim{Theorem 10}  If ${\Cal H}$ is a nonempty subset of ${\Cal
E}_\lambda$ and for each $\theta < \lambda$ $\{ h \in {\Cal H} : \text{cr} 
\,h < \theta \}$ is finite, then $\text{tp} \,(\text{cr} \,{\Cal A}_{\Cal 
H} ) = \omega$.
\endproclaim

\demo{Proof}  Since for every $\theta < \lambda$, $\text{cr} \,{\Cal
A}_{\Cal H} \cap \theta = \text{cr} \,{\Cal A}_{\{ h \in {\Cal H} :
\text{cr} \,h < \theta \}} \cap \theta$, we are reduced to the case
where ${\Cal H}$ is finite.

For $\vec p \in ({\Cal E}_\lambda )^{< \omega}$, $\vec p = \langle p_0
,p_1 ,\ldots ,p_n \rangle$, let $e(\vec p ) = p_0 p_1 \cdots p_n$.
Let ${\Cal H}^\ast = \{ e(\vec h ) : \vec h \in ({\Cal H} )^{< \omega} \}$.

\proclaim{Lemma 10.1}  $\text{tp} \,(\text{cr} \,{\Cal H}^\ast ) = \omega$
and $\sup (\text{cr} \,{\Cal H}^\ast ) = \lambda$
\endproclaim

\demo{Proof}  For each $\delta \in \text{cr} \,{\Cal H}^\ast$ pick $\vec
h_\delta \in ({\Cal H} )^{< \omega}$ of minimal length such that
$\text{cr} \,(e(\vec h_\delta )) = \delta$.  Then if $\vec h_\delta =
\langle h_0 ,h_1 ,\ldots ,h_n \rangle$ and $m < n$, then $\text{cr}
\,(h_0 h_1 \cdots h_m ) < \delta$.  Namely, if $m < n$ is such that
$\text{cr} (h_0 h_1 \cdots h_m ) \geq \delta$, then $\text{cr} \,(h_0
h_1 \cdots h_m ) = \delta$ would contradict minimality, and $\text{cr}
\,(h_0 h_1 \cdots h_m ) > \delta$ would imply $\text{cr} \,(h_{m+1}
\cdots h_n ) = \delta$, contradicting minimality.

Now suppose for some $\theta < \lambda$ that $\text{cr} \,{\Cal H}^\ast \cap
\theta$ is infinite.  Then $T$, the set of all initial segments of the
$\vec h_\delta$'s ($\delta \in \text{cr} \,{\Cal H}^\ast \cap \theta$) is a
downwards closed, finitely branching, infinite tree such that for each
$\vec h \in T$, $\text{cr} \,(e(\vec h )) < \theta$.  Then an infinite
branch through $T$ contradicts Theorem 9.
\enddemo

Let $\gamma_n^{\Cal H} = \gamma_n$ be the $n$th member of $\text{cr}
\,({\Cal A}_{\Cal H} ) = \text{cr} \,({\Cal P}_{\Cal H} )$ ($n = 0,1,2,
\ldots$).

\proclaim{Lemma 10.2}  If $\vec p \in ({\Cal P}_{\Cal H} )^{< \omega}$ 
has length
$2^n$, then for some initial segment $\vec q$ of $\vec p$, $\text{cr}
\,(e(\vec q )) \geq \gamma_n$.
\endproclaim

\demo{Proof}  The case $n=0$ is immediate.  Suppose it is true for
$n-1$.  Let $\vec p = \langle p_1 ,p_2 ,\ldots ,p_{2^n} \rangle$.  Pick
$i \leq 2^{n-1}$ such that $\ell = p_1 \cdots p_i$ has critical point
$\geq \gamma_{n-1}$; if $\text{cr} \,\ell > \gamma_{n-1}$, we are done
so suppose $\text{cr} \,\ell = \gamma_{n-1}$.  Pick $k$, $1 \leq
k\leq 2^{n-1}$ minimal so that $\text{cr} \,(p_{i+1} \cdots
p_{i+k} ) \geq \gamma_{n-1}$.  Then by Lemma 1 (i) and (iii)a, $p_1
\cdots p_{i+k} = \ell p_{i+1} \cdots p_{i+k} \overset \rho \to{=} 
\ell (p_{i+1} \cdots p_{i+k} )$ for some $\rho > \gamma_{n-1}$.  Thus
$\text{cr} \,(p_1 \cdots p_{i+k} ) > \gamma_{n-1}$, since $\text{cr}
\,(\ell (p_{i+1} \cdots p_{i+k} )) = \ell (\text{cr} \,(p_{i+1} \cdots
p_{i+k} )) > \gamma_{n-1}$ and $\rho > \gamma_{n-1}$.  So $\text{cr}
\,(p_1 \cdots p_{i+k} ) \geq \gamma_n$.
\enddemo

>From Lemma 10.1, it suffices for the theorem to show

\proclaim{Lemma 10.3} $\text{cr} \,{\Cal H}^\ast = \text{cr} \,{\Cal
P}_{\Cal H}$
\endproclaim

\demo{Proof}  Note that for $p \in {\Cal P}_{\Cal H}$, $\gamma <
\lambda$, either there is an $h \in {\Cal H}$ with $p \overset \gamma
\to{=} h$ or for some $h \in {\Cal H}$ and $q \in {\Cal P}_{\Cal H}$, 
$p \overset
\gamma \to{=} qh$.  Namely, by Lemma 7 (iv), $p \overset \gamma \to{=}
r$ for some $r \in {\Cal A}_{\Cal H}$.  Then either $r \in {\Cal H}$ or
$r$ is of the form $a_0 (a_1 (\cdots (a_n h)))$ for some $h \in {\Cal
H}$, $a_i \in {\Cal A}_{\Cal H}$; in the latter case take $q = a_0 \circ
a_1 \circ \cdots \circ a_n$.

We need to show that $\gamma_n$ is the $n$th member of $\text{cr}
\,{\Cal H}^\ast$.  Suppose this fails for some $n$, and pick $p \in
{\Cal P}_{\Cal H}$ with $\gamma_n = \text{cr} \,p$.  Fix $\gamma >
\gamma_n$, and let $p_0 = p$.  Inductively pick $p_i \in {\Cal P}_{\Cal
H}$ (for all $i \leq 2^{n+1}$ if possible) such that for each $m <
2^{n+1}$,  $p_{m+1} h_m \overset \gamma \to{=} p_m$ for some $h_m \in
{\Cal H}$.  If $m < 2^{n+1}$ and $p_0 ,\ldots ,p_m$ are chosen but
$p_{m+1}$ can't be, then $p_m \overset \gamma \to{=} h$ for some $h$,
whence $p \overset \gamma \to{=} hh_{m-1} \cdots h_1 h_0$ and $\text{cr}
\,p \in {\Cal H}^\ast$.  If the $p_i$'s can be chosen for all $i \leq
2^{n+1}$, then $p \overset \gamma \to{=} p_{2^{n+1}} h_{2^{n+1} -1}
\cdots h_1 h_0$.  By Lemma 10.2, some $\text{cr} \,(p_{2^{n+1}}
h_{2^{n+1} -1} \cdots h_i ) \geq \gamma_{n+1}$, where $i > 0$ since
$\text{cr} \,p < \gamma_{n+1}$.  But then $p \overset \gamma \to{=}
(p_{2^{n+1}} h_{2^{n+1}-1} \cdots h_i )(h_{i-1} \cdots h_1 h_0 )
\overset \gamma \to{=} h_{i-1} \cdots h_1 h_0$, giving $\gamma_n \in
\text{cr} \,{\Cal H}^\ast$.

This completes the proof of Theorem 10; note that it and the proof of
Lemma 10.3 give (for ${\Cal H}$ as hypothesized in Theorem 10) that for
all $\gamma < \lambda$ and $p \in {\Cal P}_{\Cal H}$ there's  a $q \in
{\Cal H}^\ast$ with $p \overset \gamma \to{=} q$.
\enddemo

We consider the case ${\Cal H} = \{ j\}$ below.  For $j \in {\Cal
E}_\lambda$ let $\langle \gamma_n : n < \omega \rangle$ be the
increasing enumeration of $\text{cr} \,{\Cal A}_j$; so $\gamma_0 =
\text{cr} \,j$.  Let $j_{(0)} = \text{id}$ (recall $\text{id} \notin
{\Cal E}_\lambda$), $j_{(1)} = j$ and in general $j_{(n+1)} = j_{(n)} j$.
Let $P_n = \{p \overset \ast \to\bigcap V_{\gamma_n} : p \in {\Cal P}_j
\}$.  Then from Theorem 10 and the remarks before Lemma 2, if $\cdot$
and $\circ$ are the operations on $P_n$ induced by the projection map $p
\rightarrow p \overset \ast \to\bigcap V_{\gamma_n}$, and $\pi_{mn}$ ($m
> n$) is the projection map from $P_m$ to $P_n$, then $\overline {\Cal
P}_j$, ${\Cal P}_j$, ${\Cal A}_j$ are subalgebras of the inverse limit
of $(\langle P_n : n < \omega \rangle , \langle \pi_{mn} : m > n \rangle
)$.  Parts (iii) and (iv) of the next theorem imply that $\text{Card}
\,P_n  = 2^n$.

\proclaim{Theorem 11}  (i) $j_{(2^n )} (\text{cr} \,j_{(2^n )} ) =
\text{cr} \,j_{(2^{n+1} )}$ (so in particular $\langle \text{cr}
\,j_{(2^n )} : n < \omega \rangle$ is increasing).  And $\text{cr}
\,j_{(m)} = \text{cr} \,j_{(2^n )}$ where $n$ is greatest such that $2^n$
divides $m$.
\item {(ii)} $\text{cr} \,j_{(2^n )} = \gamma_n$.
\item {(iii)} if $i, m < 2^n$ and $i \not= m$, then $j_{(i)} \overset
\ast \to\bigcap V_{\gamma_{n-1} +1} \not= j_{(m)} \overset \ast
\to\bigcap V_{\gamma_{n-1} +1}$.
\item {(iv)} $P_n = \{ j_{(i)} \overset \ast \to\bigcap V_{\gamma_n} : i
< 2^n \}$.
\endproclaim

\demo{Proof}  (i) By induction on the maximum $\ell$ of the subscripts
appearing in the expressions.  For $\ell \leq 1$ there is nothing to
prove.

For the first clause, assuming (i) holds for all $\ell < 2^{n+1}$, we
have by induction that $\text{cr} \,j_{(2^n +i)} < \text{cr} \,j_{(2^n
)}$ ($0 < i < 2^n$).  Thus $j_{(2^n )} (\text{cr} \,j_{(2^n )} ) =
j_{(2^n )} (j_{(2^n -1)} (\text{cr} \,j )) = j_{(2^n )+(2^n -1)}
(\text{cr} \,j)$ which, by Lemma 1 (iii)b, $= j_{(2^{n+1} -1)} (\text{cr}
\,j) = \text{cr} \,j_{(2^{n+1)}}$.  The second clause of (i) follows
using the induction hypothesis.

(ii) By Lemma 10.3 $\{ \gamma_n : n < \omega \} = \{ \text{cr} \,j_{(i)}
: 0 < i < \omega \}$, which equals $\{ \text{cr} \,j_{(2^n )} : n <
\omega \}$ by (i).  Since $\langle \text{cr} \,j_{(2^n )} : n < \omega
\rangle$ is increasing, $\text{cr} \,j_{(2^n )} = \gamma_n$.

(iii) By induction on $n$.  Given $i < m < 2^n$, if $m \not= 2^{n-1} + i$,
then the result follows from the induction hypothesis and the fact that
$j_{(2^{n-1} +k)} \overset \gamma_{n-1} \to{=} j_{(k)}$ 
($k < 2^{n-1}$).
If $m = 2^{n-1} +i$ and $j_{(i)} \overset \gamma_{n-1} +1 \to{=}
j_{(m)}$, then $2^{n-1} -i$ applications of the rule $p \overset
\gamma_{n-1} +1 \to{=} q \Rightarrow pj \overset \gamma_{n-1} +1 \to{=}
qj$ yield $j_{(2^{n-1} )} \overset \gamma_{n-1} +1 \to{=} j_{(2^n )}$,
then contradicting $\text{cr} \,j_{(2^{n-1} )} = \gamma_{n-1} < \text{cr}
\,j_{(2^n )}$.

(iv)  Given $p \in {\Cal P}$.  By Lemma 10.3 there is an $i$ with $p
\overset \gamma_n \to{=} j_{(i)}$; pick $i$ minimal.  If $\text{cr} \,p
\geq \gamma_n$, $i = 0$.  If $\text{cr} \,p < \gamma_n$, then if $i \geq
2^n$, we would have by Lemma 10.2 an $m < i$ with $\text{cr} \,j_{(m)}
\geq \gamma_n$, whence $p \overset \gamma_n \to{=} j_{(i-m)}$.  So $i <
2^n$.

Recall that $A$ is the set of terms $w$ in the operation $\cdot$, in one
variable $x$, and for $j \in {\Cal E}_\lambda$, $w[j]$ is the member of
${\Cal A}_j$ induced by the homomorphism which sends $[x]$ to $j$.

\proclaim{Theorem 12}  $\{ \langle u,w\rangle \in A \times A : \text{cr}
\,(u[j]) < \text{cr} \,(w[j])\}$ is recursive and independent of $j$.
\endproclaim

\demo{Proof}  The tables of the finite algebras $P_k$ can be defined
without reference to large cardinals.  Define a binary operation
$\ast_k$ on $2^k = \{ 0,1,\ldots ,2^k -1\}$ by
$$\aligned m \ast_k 0 &= 0\\
m \ast_k 1 &= m+1 ~(\text{mod} \,2^k )\endaligned$$
$$m \ast_k i = [m \ast_k (i-1)] 
\ast_k [m \ast_k 1].\quad (1 < i < 2^k )$$

To see that this constitutes a recursive definition of an operation, note that 
the equations
imply $0 \ast_k n = n$, $2^k -1\ast_k n = 0$, and (by induction on $m = 0$, 
$2^k -1 ,2^k -2,\ldots ,1$)
that $m \ast_k n = 0$ or $m \ast_k n > m$.  By the latter clause and
a similar induction, $m \ast_k n$ is (uniquely) determined for all
$m, n < 2^k$.

We claim that for $m, n, \ell < 2^k$, $m \ast_k n = \ell$ 
if and only if $j_{(m)} j_{(k)}
\overset \gamma_k \to{=} j_{(\ell )}$, i.e., $(2^k ,\ast_k )
\cong (J_k ,\cdot )$.  Namely, $(J_k ,\cdot )$ satisfies the 
first and third
equations since $j_{(0)} = \text{id}$ and $(J_k ,\cdot )$ is left
distributive, and the second holds by definition unless $m = 2^k -1$, 
in which case $j_{(2^k -1)} j = j_{(2^k )} \overset \gamma_k
\to{=} j_{(0)}$.
If $v \in A$, $k < \omega$, the unique $i < 2^k$ such that $v[j]
\overset \gamma_k \to{=} j_{(i)}$ may thus be computed in the algebra
$(2^k ,\ast_k )$.  And $\text{cr} \,(v[j]) \geq \gamma_k$ 
iff $i = 0$, and if $0 < i < 2^k$, then $\text{cr} \,(v[j]) =
\gamma_\ell$, where $2^\ell$ is the largest power of 2 dividing $i$.
For sufficiently large $k$ the number $i$ thus associated to $w$ in $(2^k
,\ast_k )$ is nonzero (using the assumption $\exists \,\lambda {\Cal
E}_\lambda \not= \emptyset$ and its consequence $\sup \{ \text{cr} \,j_{2^n} :
n < \omega \} = \lambda$).  This gives an algorithm for comparing
$\text{cr} \,w[j]$ with $\text{cr} \,u[j]$.
\enddemo

A number of statements about ${\Cal A}_j$ and its members can be coded
up as staatements about the $(2^k ,\ast_k )$'s, thus reducing them to
problems in finite combinatorics, and raising the question whether
those versions are provable in ZFC.  Wehrung ([15]) proved (in ZFC) that
for $m$ a positive integer, the free left-distributive algebra on one
generator $x$ satiisfying $x_{(m+1)} = x$ is $(2^k ,\ast_k )$, when $k$
is greatest such that $2^k$ divides $m$.  Regarding the result that
$\text{tp} (\text{cr} {\Cal A}_j ) = \omega$, it can be phrased finite
combinatorially as the statement that each term $w$ in $A$, evaluated in
each $(2^k ,\ast_k )$, is nonzero for eventually all
$k$, or as the statement that for all $i \geq 1$  there's a $k$ with $1
\ast_k i \not= 0$.  These (and other) versions aren't known to be
provable in ZFC.

A result of Woodin, using upward Easton forcing, is that it's
consistent (relative to the consistency that some ${\Cal E}_\lambda
\not= \emptyset$) that there are $k, \ell \in {\Cal E}_\lambda$ with $k
\not= \ell$ but $k \restriction \lambda = \ell \restriction \lambda$.
However, if for some $j \in {\Cal E}_\lambda$, $k, \ell \in \overline
{\Cal A}_j$, then if $k \restriction \lambda = \ell \restriction
\lambda$, then $k = \ell$.  Namely, let $\text{cr} (j)  = \kappa$,
$<_\kappa$ a well ordering of $V_\kappa$ of length $\kappa$.  Then
$<_\lambda = \underset n < \omega \to\bigcup j^n (<_\kappa )$, the
$\omega \text{th}$ iterate of $<_\kappa$, is a well ordering of
$V_\lambda$ of length $\lambda_j$ and by induction on ${\Cal A}_j$ one
has that for each $p \in  \overline {\Cal A}_j$, $p (<_\lambda ) =
<_\lambda$.  Then if $k(y) \not= \ell (y)$ and $y$ is the $\alpha
\text{th}$ member of $V_\lambda$ under $<_\lambda$, then $k(\alpha )
\not= \ell (\alpha )$, as desired.

Using Theorem 11, for $j \in {\Cal E}_\lambda$ embeddings $p, q \in
\overline {\Cal A}_j$ can be constructed with $p \restriction \text{cr}
({\Cal A}_j ) =  q \restriction \text{cr} ({\Cal A}_j )$ and $p \not=
q$.  Members of ${\Cal A}_j$, though, are determined by their
restrictions to $\text{cr} ({\Cal A}_j )$.

\proclaim{Theorem 13} If $j \in {\Cal E}_\lambda$, $p,q \in {\Cal A}_j$,
$p \not= q$, then $p \restriction \text{cr} \,{\Cal A}_j \not= q
\restriction \text{cr} \,{\Cal A}_j$.
\endproclaim

We remark that counterexamples show that ``$\text{cr} \,{\Cal A}_j$''
can't be replaced in the theorem by ``$\text{cr} \,{\Cal A}_n \cap
\gamma_n$, where $\gamma_n$ is least with $p \overset \gamma_n  
\to{\not=} q$,'' nor can it be
replaced by ``$\{ j^n (\text{cr} \,j) : n < \omega \}$.''

\demo{Proof of the Theorem}
\enddemo

\proclaim{Lemma 13.1}  If $k \in {\Cal E}_\lambda$, $\rho < \lambda$,
then for some $\delta < \lambda$ with $\text{cr} \,\delta > \rho$ there
is a $\theta \in \text{cr} \,{\Cal A}_k$ with $k(< \delta ) < \theta <
k(\delta )$.
\endproclaim

\demo{Proof}  Let $k^{(1)} = kk$, $I_n = I_n (k^{(1)} ,k)$ ($n <
\omega$).  Then $\underset n \to\sup (\text{cr} \,I_{2n+1} ) = \lambda$ 
by Lemma 7 (iii).  Pick an $n$ such that $\mu = \text{cr} \,I_{2n+1} >
\rho$.  Let $\delta = k(< \mu )$, of cofinality $\mu > \rho$.
\enddemo

We have $k (< \delta ) = I_{2n+2} (< \mu )$.  Namely, $k (< \delta ) = k
(< k(< \mu )) = k^{(1)} (< k (< \mu )) = (I_1 \circ I_0 ) (< \mu ) =
(I_{2n+2} \circ I_{2n+1} )(< \mu ) = I_{2n+2} (< I_{2n+1} (< \,\mu )) =
I_{2n+2} (< \mu )$ since $\mu = \text{cr} \,I_{2n+1}$.

We have $k (\delta ) = [I_1 (< I_0 \mu )] > [I_2 (< I_1 \mu )] > [I_3 (<
I_2 \mu )] > \cdots > [I_{2n+2} (< I_{2n+1} \mu )]$.  Namely $k (\delta
) = k (k (< \mu )) = k^{(1)} (< k \mu ) = I_1 (< I_0 \mu )$.  To see
that $[I_k (< I_{k -1} \mu )] \gg [I_{k+1} (< I_k \mu
)]$ ($k\leq 2n+1$), we have that $I_{k-1} (\mu ) > I_{k
-1} (< \mu )$ since $\text{cr} \,I_{k-1} < \mu$ and $\mu$ is
regular, so $[I_k (< I_{k-1} \mu )] > I_k (I_{k-1}
(<\mu )) = I_{k+1} (< I_k \mu )$.

So $k(< \delta ) = I_{2n+2} (< \,\mu ) < I_{2n+2} (\mu )$ (since
$\text{cr} \,I_{2n+2} = \text{cr} \,k < \mu ) < [I_{2n+2} (< I_{2n+1}
\mu )]$ (since $\text{cr} \,I_{2n+1} = \mu ) < [I_{2n+1} (< I_{2n} \mu
)] < \cdots < k (\delta )$.  Therefore $\theta = I_{2n+1} (\mu )$ is as
desired.

For $w \in A$, let $\bar w = w[j]$.

\proclaim{Lemma 13.2}  Suppose $F \subseteq A$ is finite and closed
under subterms, $n \geq \text{Card} \,F$, $\delta < \lambda$, $\text{cr}
\,\delta > \text{cr} \,\bar w$ for all $\bar w \in F$.  Then, letting
$\beta = j^n (< \delta )$, we have
\item {(i)} $\bar u (< \beta ) = \bar v (< \beta )$ (all $u, v \in F$)
\item {(ii)} $\bar u (\beta ) > \overline {uv} (\beta )$ (all $uv \in
F$)
\endproclaim

\demo{Proof}  By induction on $\text{Card} \,F$.  If $F = \{ x\}$, it is
trivial.  Suppose $F = F' \cup \{ uv\}$, where $uv \notin F'$ and $F'$
is closed under subterms.  Given $n \geq \text{Card} \,F$, $\text{cr}
\,\delta \geq \text{cr} \,\bar w$ (all $w \in F$), $\beta = j^n (<
\delta )$.

For (i) we have $\bar w (< \beta ) = \bar t (< \beta )$, all $\bar x ,t
\in F'$.  And $\overline {uv} (< \beta ) = \overline {uv} (< j^n  (< 
\delta )) =
\overline {uv} (< j (< j^{n-1} (< \delta ))) = \overline {uv} (< \bar u
(< j^{n-1} (< \delta )))$ (induction hypothesis) $= \bar u (< \bar v (<
j^{n-1} (< \delta ))) = \bar u (< j^n (< \delta )))$ (induction
hypothesis) $= \bar u (< \beta )$, as desired.

For (ii), $\bar u (\beta ) = \bar u (j^n (< \delta )) = \bar u (j(<
j^{n-1} (< \delta ))) = \bar u (\bar v (< j^{n-1} (< \delta ))) =
\overline {uv} (< \bar u (j^{n-1} (< \delta ))) > \overline {uv} (\bar u
(< j^{n-1} (< \delta )))$ (since $\text{cr} \,j^{n-1} (< \delta ) =
\text{cr} \,\delta \geq \text{cr} \,\bar u ) = \overline {uv} (\beta )$.
\enddemo

To prove the theorem, let $p, q \in {\Cal A}_j$, $p \not= q$.  Then by
Theorem 5 (i), we have, say, $p <_L q$.  Take $w \in A$, $\bar w = p$,
$w = u_0 u_1 \cdots u_m$, $\bar u_0 = q$.  Let $F$ be the set of all
subterms of $w$.  Applying Lemma 13.1 to $k = p$, there is a $\delta <
\lambda$ with $\text{cr} \,\delta \geq \sup \{ \text{cr} \,\bar w : w
\in F\}$ and a $\theta \in \text{cr} \,{\Cal A}_p$ (whence in $\text{cr}
\,{\Cal A}_j$) with $p (< \delta ) < \theta < p(\delta )$.

Let $n = \text{Card} \,F$, $\beta = j^n (< \delta )$.  By Lemma 13.2,
$q(\beta ) = \bar u_0 > \overline{u_0 u_1} (\beta ) > \cdots > \overline
{u_0 u_1 \cdots u_m} (\beta ) = p(\beta )$, and $\bar u (< \beta ) =
\bar v (< \beta )$ (all $u, v \in F$).  Claim there is a $\gamma \in
\text{cr} \,{\Cal A}_j$ with $p (< \beta ) < \gamma < p(\beta )$.
Namely, since $p (< \delta ) < \theta < p (\delta )$, 
we have $p(< \beta ) = p(< j^n (< \delta )) = pj^n (< p (<
\delta )) < pj^n (\theta ) < pj^n (< p\delta ) = p(j^n (< \delta )) =
p(\beta )$.  Take $\gamma = pj^n (\theta )$.

We claim that for $i < m$, $\overline {u_0 \cdots u_i} (\gamma ) >
\overline {u_0 \cdots u_{i+1}} (\gamma )$.  Namely $\overline {u_0
\cdots u_i} (\gamma ) > \overline {u_0 \cdots u_i} (p(< \beta )) =
(\overline {u_0 \cdots u_i} )(\overline {u_{i+1}} (< \beta )) =
\overline {u_0 \cdots u_i u_{i+1}} ( < \overline {u_0 \cdots u_i} (\beta
)) > \overline {u_0 \cdots u_{i+1}} (p(\beta )) > \overline {u_0 \cdots
u_{i+1}} (\gamma )$.  Thus $q(\gamma )= \bar u_0 (\gamma ) > \overline
{u_0 \cdots u_m} (\gamma ) = p(\gamma )$, proving the theorem.
\vfill\eject

\baselineskip=12pt

\centerline {\bf References}
\vskip 1em
\item {[1]} P. Dehornoy, Free distributive groupoids, {\it Journal of
Pure and Applied Algebra} {]bf 61} (1989), 123--146.
\vskip 10pt
\item {[2]} P. Dehornoy, Sur la structure des gerbes libres, {\it CRAS
Paris} {\bf 309-I} (1989), 143--148.
\vskip 10pt
\item {[3]} P. Dehornoy, Braid groups and left distributive operations,
preprint.
\vskip 10pt
\item {[4]} R. Dougherty, Critical points of elementary embeddings,
handwritten notes, 1988.
\vskip 10pt
\item {[5]} A. Kanamori, W. Reinhardt and R. Solovay, Strong axioms of
infinity and elementary embeddings, {\it Ann. Math. Logic} {\bf 13}
(1978), 73--116.
\vskip 10pt
\item {[6]} A. Kanamori, Large cardinals, Springer.
\vskip 10pt
\item {[7]} K. Kunen, Elementary embeddings and infinitary
combinatorics, {\it J. Symb. Logic} {\bf 36} (1971), 407--413.
\vskip 10pt
\item {[8]} D. Larue, On braid words and irreflexivity, {\it Algebra
Universalis}
\vskip 10pt
\item {[9]} R. Laver, The left distributive law and the freeness of an
algebra of elementary embeddings, {\it Advances in Mathematics} {\bf 91}
(1992), 209--231.
\vskip 10pt
\item {[10]} R. Laver, A division algorithm for the free left
distributive algebra, {\it Proc. Helsinki 1990 ASL Meeting},
\vskip 10pt
\item {[11]} D. A. Martin and J. Steel, A proof of projective
determinacy, {\it Jour. Amer. Math. Soc}. {\bf 2} (1989), 71--125.
\vskip 10pt
\item {[12]} D. A. Martin and J. Steel, Iteration trees, {\it Jour.
Amer. Math. Soc}.
\vskip 10pt
\item {[13]} W. Mitchell, Sets constructible from sequences of
ultrafilters, {\it JSL} {\bf 39} (1974), 57--66.
\vskip 10pt
\item {[14]} J. Steel, The well foundedness of the Mitchell order,
preprint.
\vskip 10pt
\item {[15]} F. Wehrung, Gerbes primitives, {\it CRAS Paris}, 1991.
\vskip 2em
\noindent University of Colorado, Boulder.

\end